# In Polytopes, Small Balls about Some Vertex Minimize Perimeter

Frank Morgan
Department of Mathematics and Statistics
Williams College
Williamstown, MA 01267
Frank.Morgan@williams.edu

## Abstract

In (the surface of) a convex polytope $P^n$ in $\mathbf{R}^{n+1}$, for small prescribed volume, geodesic balls about some vertex minimize perimeter. This revision corrects a mistake in the mass bound argument in the proof of Theorem 3.8.

# 1. Introduction

Even in smooth Riemannian manifolds, there are relatively few examples of explicitly known regions which minimize perimeter for prescribed volume (see [M1, 13.2], [HHM]). One general result is that for small volume, a perimeter-minimizing region is a nearly round ball where the scalar curvature is large ([K], [MJ, Thm. 2.2], or in 3D [Ros, Thm. 18], with [Dr]).

For singular ambients, one result in general dimensions is that in smooth cones of positive Ricci curvature, geodesic balls about the apex are perimeter minimizing [MR, Cor. 3.9]. Examples of Cotton *et al.* ([GNY], [CFG]) include the surface of the cube in $\mathbf{R}^3$, with perimeter-minimizing regions as pictured in Figure 1. In particular, for small prescribed area, geodesic balls about a vertex minimize perimeter. The analogous result in higher dimensions has remained open (see [MR, Rmk. 2.3]). Our Theorem 3.8 proves that

> *In a polytope $P^n$ in $\mathbf{R}^{n+1}$, for small prescribed volume, geodesic balls about some vertex minimize perimeter.*

(Note that by "polytope" we refer to the *boundary* of the compact, convex solid body.) For $n \geq 3$, $P^n$ has singularities of positive dimension, and the perimeter-minimizing spheres pass through these singularities.

**1.1. The proof.** The proof has the following elements. Let C be the tangent cone at a vertex, with link a spherical polytope $Q^{n-1}$ in $\mathbf{S}^n$.

(1)     By the Levy-Gromov isoperimetric inequality, the isoperimetric profile for Q is bounded below by the isoperimetric profile for $\mathbf{S}^{n-1}$. This step requires smoothing (Lemma 2.1) and approximation.

(2)     By an isoperimetric comparison theorem 3.2 for products and cones (Ros after Barthe and Maurey), balls about the apex of C are perimeter minimizing (Theorem 3.3).

(3)     For perimeter-minimizing surfaces in $P^n$ with prescribed volumes approaching zero, strong control of mean curvature after Almgren [Alm] (Lemma 3.5) is obtained by considering a limit of renormalizations to unit volume. A volume concentration argument shows that the limit is nonzero.

(4)     It follows from (3) that the sum of the diameters of components of perimeter-minimizing regions is small (Lemma 3.6). The proof involves a generalization of lower mass bounds for minimizers of elliptic integrands to allow volume constraints. As a result we may assume that components lie in a neighborhood of the vertices.

(5)     By (2), we may assume that they consist of balls about vertices. A single ball is best.

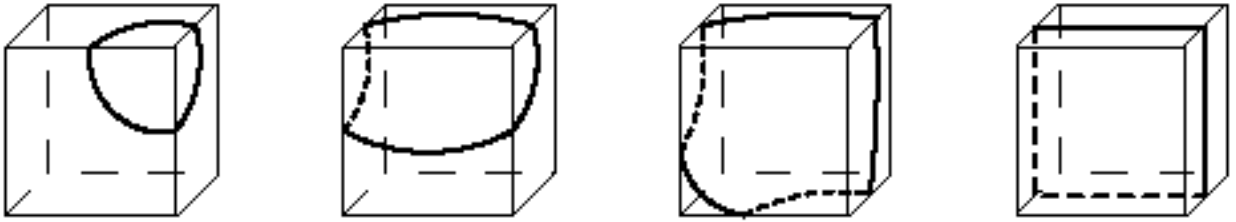

Figure 1
Perimeter-minimizing regions in the surface of the cube [CFG, Fig. 1].

**1.2. Solid polytopes.** Our results also hold in the less singular category of *solid* polytopes (3.11).

**1.3. Acknowledgments.** I learned about the isoperimetric comparison theorem for products 3.2 in Ros's 2001 Clay-MSRI lectures [Ros], and profited from conversations with Ros, Ritoré, McCuan, and Bayle while visiting Granada in May and June, 2003. I am grateful for proddings for clarity from the editor and referees, which caused me to find and correct a mistake in the mass bound argument in the proof of Theorem 3.8 with new ideas for generalizations. This research was partially supported by a National Science Foundation grant.

## 2. Smoothing

The following lemma will be used in the proof of Theorem 3.3.

**2.1. Smoothing Lemma.** *Given a (convex) polytope $Q^{n-1}$ in $S^n$, there are smooth convex hypersurfaces $Q_\varepsilon$ which converge smoothly to Q off the singular set with volumes converging to that of Q.*

*Proof.* The corresponding result in $R^n$ is standard and easy (see e.g. [Gh, Note 1.3]). Indeed, fix an origin in $R^n$ inside Q, take the convex function F which is the fraction of the distance from the origin to Q along a straight line, let $F_\varepsilon$ be a smoothing of F by symmetric convolution, and let $Q_\varepsilon = \{F_\varepsilon = 1\}$. For convolution, start with a nice smooth symmetric function $\varphi$ supported in (0,1) with integral 1 and let $\varphi_\varepsilon(x) = \varphi(x/\varepsilon)/\varepsilon$.

In $S^n$, we may assume that Q lies in the northern hemisphere H, where distance from the north pole is a convex function. Take the convex function F on H which is the fraction of the distance from the north pole to Q along a geodesic. Choose $\zeta \leq \varepsilon$ such that on a $\zeta$ neighborhood of $\{F = 1-\varepsilon\}$, $1-2\varepsilon \leq F \leq 1-\varepsilon$. Let $F_\varepsilon$ be the convolution of F with $\varphi_\zeta$, which is convex. Let $Q_\varepsilon = \{F_\varepsilon = 1-2\varepsilon\}$.

## 3. Isoperimetric Regions

**3.1. Definitions.** The cone C over a subset M of the unit sphere in Euclidean space is defined as

$$C = \{tx : x \in M, t \geq 0\}.$$

For a general smooth Riemannian manifold M, the cone (minus the negligible apex) may be defined as the warped product $(0,\infty) \times_f M$ with $f(t) = t$, i.e., as $(0,\infty) \times M$ with the metric

$$ds^2 = dt^2 + t^2\, dm^2.$$

A *density* on a Riemannian manifold is a positive function used to weight volume and perimeter (see [M2]). Although we will use only constant densities, Theorem 3.2 holds for variable smooth densities. For the more general context of "mm spaces" (metric spaces with measures), see Gromov [Gr2].

The *isoperimetric profile* I(M)(V) of a space M gives the minimum (or infimum) perimeter required to enclose volume V.

The following useful theorem of Ros says that isoperimetric comparison of a manifold $M_2$ with the sphere is preserved under product with another manifold $M_1$. It follows earlier comparisons with Gauss space by Barthe and Maurey. The idea is Schwarz symmetrization: compare a region in a product $M_1 \times M_2$ with the region in $M_1 \times S^n$ obtained by replacing vertical slices (contained in $M_2$) with (n–1)-balls in $S^n$.

**3.2. Isoperimetric Comparison Theorem for Warped Products** [Ros, Thm. 3.7]. *Let $M_1$, $M_2$ be smooth Riemannian manifolds with densities, with $|M_2| = |S^n|$. Consider smooth warped products $M_1 \times_f S^n$, $M_1 \times_f M_2$. If $I(M_2) \geq I(S^n)$, then $I(M_1 \times_f M_2) \geq I(M_1 \times_f S^n)$. In particular, if for some region $R \subset M_1$ the product $R \times S^n$ is minimizing, so is $R \times M_2$, uniquely if $I(M_2) > I(S^n)$ except at the endpoints.*

*Remarks.* Although Ros does not consider *warped* products, the proof is the same (and the uniqueness statement follows immediately from the proof).

Theorem 3.2 provides an immediate proof of a theorem of Morgan and Ritoré. Our generalization in part B requires our Smoothing Lemma 2.1.

**3.3. Theorem. A** [MR, Cor. 3.9]. *For $n \geq 3$, let C be a cone over a smooth, closed, connected Riemannian manifold $M^{n-1}$ of Ricci curvature at least $n-2$. Then geodesic balls about the apex are perimeter minimizing for given volume, uniquely unless $C = \mathbf{R}^n$.*

**B.** *Likewise, for $n \geq 2$, let C be the cone over any (convex) spherical geodesic polytope $Q^{n-1} \subset S^n$. Then geodesic balls about the apex are perimeter minimizing.*

*Remark.* A similar, trivial result is that in $\mathbf{R}^n$ modulo a finite group $\Gamma$ of isometries leaving a point $p$ fixed, geodesic balls about $p$ are perimeter minimizing, equivalent to the obvious statement that among $\Gamma$-invariant objects in $\mathbf{R}^n$, geodesic balls about $p$ are perimeter minimizing. This observation applies to any perimeter-minimizing geodesic ball about a point $p$ in any Riemannian manifold.

*Proof of Theorem 3.3.* Renormalize M with constant density to make $|M| = |S^{n-1}|$. By the curvature hypothesis, the Levy-Gromov isoperimetric inequality ([Gr1, 2.2], [BZ, 34.3.2], or [Ros, Sect. 2.5]) says that the isoperimetric profiles satisfy

(1) $$I(M) \geq I(S^{n-1}).$$

Note that C is the warped product $(0,\infty) \times_f M$ with $f(t) = t$, and that $(0,\infty) \times_f S^{n-1}$ is just $\mathbf{R}^n$. Hence by Theorem 3.2, since balls about the origin are perimeter minimizing in $\mathbf{R}^n$, balls about the apex are perimeter minimizing in C. If there is another minimizer, then equality holds nontrivially in (1), and $M = S^{n-1}$, i.e., $C = \mathbf{R}^n$.

Likewise given any (convex) spherical geodesic polytope $Q^{n-1} \subset S^n$ with $n \geq 3$, given a volume V, use Lemma 2.1 to obtain a smoothing $Q_\varepsilon$ with sectional curvature at least 1 and hence Ricci curvature at least $n-2$ (and normalized by a constant density to have the same volume as $S^{n-1}$). Then as before

(2) $$I(C_\varepsilon) \geq I(\mathbf{R}^n).$$

Hence

$$I(C)(V) \geq \limsup I(C_\varepsilon)(V) \geq \liminf I(C_\varepsilon)(V) \geq I(\mathbf{R}^n)(V),$$

and balls around the apex are perimeter minimizing.

For the easier case n = 2, when all cones are isometric to circular cones, geodesic balls about the apex are well known to be uniquely perimeter minimizing. See [HHM, Sect. 8]; alternatively, apply Theorem 3.2 after noting that since $|Q| \leq |\mathbf{S}^1|$, after renormalization $I(Q) \geq I(\mathbf{S}^1)$.

The following Lemma 3.4 is useful in obtaining uniform estimates in the proof of the main Theorem 3.8.

**3.4. Lemma.** *Fix a simplex $\Delta$ in $\mathbf{R}^n$. For each integer $N \geq 2$, slice $\Delta$ by N equally spaced hyperplanes parallel to each face, the first containing the face, the last through the opposite vertex. Among the many pieces of $\Delta$ which arise under such slicing for all such N, at most n shapes occur (up to translation and homothety).*

*Remark.* In $\mathbf{R}^2$, slicing a triangle yields two shapes: the triangle and its central inversion (an upside-down triangle). In $\mathbf{R}^3$, slicing a tetrahedron yields three shapes: the tetrahedron, its inversion, and an octahedron (fully truncated tetrahedron).

*Proof of Lemma 3.4.* By a linear transformation, we may assume that $\Delta$ is regular, with unit normals $a_1 + \ldots + a_{n+1} = 0$, normalized so that

(1) $\qquad\qquad a_i \bullet a_j = -1$, except that $a_i \bullet a_i = n$.

Such a regular simplex $\Delta$ is given by

$$\Delta = \{ -1 < a_i \bullet x < 0, \text{ except } 0 < a_{n+1} \bullet x < 1 \}.$$

Our slicing yields pieces of the form

$$\{ -\frac{k_i}{N} < a_i \bullet x < -\frac{k_i-1}{N}, \text{ except } \frac{k_{n+1}-1}{N} < a_{n+1} \bullet x < \frac{k_{n+1}}{N} : 1 \leq k_i \leq N \},$$

or by scaling and replacing $k_{n+1}$ by $-k_{n+1}+1$, more generally

(2) $\qquad\qquad S_k = \{-k_i < a_i \bullet x < -k_i+1 : k_i \in \mathbf{Z}\}.$

By (1), translating $S_k$ by $a_1$ just adds $(-n, 1, 1, \ldots, 1)$ to k, and similarly for each $a_i$. Hence translations generate $\{\Sigma k_i = 0\}$, and the shape of $S_k$ is characterized by $\Sigma k_i$. By summing (2) over i, we see that for nonempty $S_k$,

$$0 < \Sigma k_i < n+1 ,$$

yielding at most n different shapes, as asserted.

**3.5. Lemma**. *Locally, let $S_\alpha$ be a weakly convergent sequence of mD isoperimetric hypersurfaces without boundary, with a nonzero limit. Then there are constants $C_1, \delta > 0$ such that for $\alpha$ large, in any $\delta$-ball B, $S_\alpha$ minimizes area plus $C_1|V|$ in comparison with other surfaces in B with the same boundary enclosing net volume V with $S_\alpha$.*

*Remarks and definitions.* This lemma holds in the general context of rectifiable currents [M1], with weak convergence in the flat norm or as measures. Isoperimetric surfaces satisfy an equilibrium condition of "weakly constant mean curvature" H, which means that for any smooth family of diffeomorphisms for which the initial rate of change of enclosed volume $dV/dt|_0$ is nonzero, the derivative of area A with respect to volume satisfies

$$\frac{dA}{dV}\bigg|_0 = mH.$$

The conclusion of the lemma is a kind of uniform mean curvature bound. The limit ambient may have singularities, as long as it is a smooth manifold at some point of S.

*Proof of Lemma 3.5.* Choose $\delta_1 > 0$ small enough so that there are two smooth $\delta_1$-balls about points of the limit S such that any third $\delta_1$-ball B is disjoint from one of them, say $B_0$. Let $\Phi_t$ be a smooth family of diffeomorphisms supported in the interior of $B_0$ such that for S, $dV/dt|_0 = 1$; let $A_0 = dA/dt|_0$. Choose $t_0 > 0$ so small that for all $-t_0 \leq t \leq t_0$, $.9 \leq dV/dt \leq 1.1$ and $|dA/dt| \leq |A_0| + 1$. Consequently, we may assume that for all $S_\alpha$, for all $-.9t_0 \leq V \leq .9t_0$, such diffeomorphisms of $S_\alpha$ alter enclosed volume with $S_\alpha$ by V and increase area at most $|V|(|A_0| + 1)/.9 = C_1 V$. Choose $\delta < \delta_1$ such that a $\delta$-ball has volume less than $.9t_0$.

Now let T be any surface which differs from $S_\alpha$ only in a $\delta$-ball B. Obtain T' by altering T in $B_0$ so that it bounds net volume 0 with $S_\alpha$, such that its area satisfies

$$|T'| \leq |T| + C_2|V|.$$

Since $S_\alpha$ is isoperimetric,

$$|S_\alpha| \leq |T'| \leq |T| + C_2|V|,$$

as desired.

The following lemma is a modification of a standard "monotonicity" argument (see [M1, Chapt. 9]).

**3.6. Lemma.** *Let C be an nD (convex) polytopal cone in $\mathbf{R}^{n+1}$. Given $C_1 > 1$, there is an $a > 0$ such that the following holds. Let S be an (n−1)-dimensional surface (rectifiable current) in a closed ball B of radius 1 about one of its points, with boundary contained in $\partial B$. Suppose that S minimizes area plus $C_1|V|$ in comparison with other surfaces in B with the same boundary enclosing net volume V with S. Then the area of S is at least a.*

*Proof.* Choose a vector in the interior of the cone. Projection f onto the orthogonal hyperplane $\Pi$ is a bijection of C with $\Pi$, which for some $C_2 > 0$, shrinks distance, area, and volume by a factor of at most $C_2$. In particular, there is a $\delta > 0$ such that the projection of a unit ball about a point in C contains a $\delta$-ball about the projection of the point.

Given S and a point in S, let T and p denote their images under f. For $r \leq \delta$, let $T_r = T \cap B(p, r)$, let $S_r$ denote the preimage of $T_r$, and let g(r) denote the area of $T_r$. For almost all r, the boundary of $T_r$ has area at most g'(r) (see [M1, Chapt. 9). By the isoperimetric theorem for the ball, this same boundary bounds a surface T' of area

$$|T'| \leq \alpha g'(r)^{(n-1)/(n-2)},$$

for some isoperimetric constant $\alpha$ depending only on dimension; together with $T_r$, T' bounds a region of mass at most

$$\beta_0[g(r) + \alpha g'(r)^{(n-1)/(n-2)}]^{n/(n-1)} \leq \beta g(r)^{n/(n-1)} + \beta g'(r)^{n/(n-2)}$$

for some isoperimetric constant $\beta_0$ and constant $\beta$. The preimage $S' = f^{-1}(T')$ has area

$$|S'| \leq C_2 |T'| \leq C_2 \alpha g'(r)^{(n-1)/(n-2)};$$

together with $S_r$, S' bounds a region of mass at most $C_2\beta g(r)^{n/(n-1)} + C_2\beta g'(r)^{n/(n-2)}$. By the minimizing property of S,

$$|S_r| \leq |S'| + C_1 C_2 \beta g(r)^{n/(n-1)} + C_1 C_2 \beta g'(r)^{n/(n-2)}$$

$$\leq C_2 \alpha g'(r)^{(n-1)/(n-2)} + C_1 C_2 \beta g(r)^{n/(n-1)} + C_1 C_2 \beta g'(r)^{n/(n-2)}.$$

Since $g(r) \leq |S_r|$, for some $C_3 > 0$,

$$g(r)\left(1 - C_3 g(r)^{1/(n-1)}\right) \leq C_3[g'(r)^{(n-1)/(n-2)} + g'(r)^{n/(n-2)}].$$

We may assume that the coefficient of g(r) for any $r \leq 1/2C_2 < 1$ is at least 1/2, since otherwise we have a lower bound as desired for g(r). Let

$$A = \{0 < r < 1 : g'(r) \leq 1\}.$$

We may assume that $|A| \geq 1/2$, since otherwise we immediately obtain a lower bounded as desired for g(r). On A,

$$g(r) \leq C_4 \, g'(r)^{(n-1)/(n-2)}.$$

Since g is monotonically increasing, integration yields

$$g(r) \geq C_5 \, r^{n-1},$$

again yielding the desired lower bound for $g(1/2C_1)$ and hence for |S|.

**3.7. Corollary.** *Given a polytope $P^n$ in $\mathbf{R}^{n+1}$ and $C_1$, $\delta > 0$, there is an $a > 0$ such that the following holds. Let S be an (n−1)-dimensional surface in P which in any $\delta$-ball B minimizes area plus $C_1|V|$ in comparison with other surfaces in B with the same boundary enclosing net volume V with S. If the unit ball B(p, 1) about a point p in S lies in the cone at a vertex of P, then the area of $S \cap B(p, 1)$ is at least a.*

*Proof.* By scaling we may assume that $\delta = 1$. Lemma 3.6 applies to the cone over each vertex of P. Let a be the minimum of the corresponding a's.

Gnepp, Ng, and Yoder ([GNY], [CFG]) proved that in the surface of a cube or regular tetrahedron, for small prescribed area, balls around a vertex are perimeter minimizing. The same interesting question in higher dimensions has remained open (see [MR, Rmk. 2.3]). The following theorem answers the question affirmatively. (Recall that by "polytope" we refer to the *boundary* of the compact, convex solid body.)

**3.8. Theorem.** *For $n \geq 1$, let $P^n$ be a polytope in $\mathbf{R}^{n+1}$. For small prescribed volume V, geodesic balls about some vertex are perimeter minimizing.*

*Remark.* Except for polyhedra in $\mathbf{R}^3$ (n = 2), our proof does not provide uniqueness, because of the use of smoothing and approximation. Of course uniqueness fails for polygons (n = 1).

*Proof.* Let R be a perimeter-minimizing region in P of small volume V. We will need some concentration of volume. Partition P into simplices and slice each simplex into small polytopes $K_i$ by N equally spaced hyperplanes parallel to each face. By Lemma 3.4, there are only finitely many shapes for the $K_i$ (independent of N, up to homothety) and for some positive constant $\beta$ (independent of N) we may choose N such that

$$\beta V \leq |K_i| \leq V.$$

We claim that *there is a constant $\delta > 0$ depending only on P such that some K satisfies*

(1) $$|R \cap K| \geq \delta V.$$

To prove this claim, first note that the least area A for small volume V satisfies $A \leq CV^{(n-1)/n}$. We may assume that $\delta < 1/2$ and $|R \cap K_i| \leq |K_i|/2$. Because there are only finitely many shapes for the $K_i$, by the relative isoperimetric inequality [M1, 12.3(1)], there is a single constant γ such that

$$|\partial R \cap K_i| \geq \gamma \, |R \cap K_i|^{(n-1)/n} \geq \gamma \, |R \cap K_i| \, / \, \max|R \cap K_i|^{1/n}.$$

Summing over i yields

$$A \geq \gamma \, V \, / \max|R \cap K_j|^{1/n},$$

$$\max|R \cap K_j|^{1/n} \geq \gamma \, V/A \geq \gamma \, V^{1/n} \, / \, C.$$

Hence some K satisfies

$$|R \cap K| \geq \delta V$$

with $\delta = \gamma^n / C^n$, proving the claim (1).

To complete the proof of the theorem, let $R_\alpha$ be a sequence of perimeter-minimizing regions with $V_\alpha$ small and approaching 0. For each $R_\alpha$ use the claim to choose $K_\alpha$ with $\beta V_\alpha \leq |K_\alpha| \leq V_\alpha$ and

$$|R_\alpha \cap K_\alpha| \geq \delta V_\alpha \geq \delta \, |K_\alpha|.$$

Rescale up to $1 = V_\alpha \geq |K_\alpha| \geq \beta$ with $K_\alpha$ centered at the origin. By compactness [M1, 9.1], we may assume that the $R_\alpha$ converge to a perimeter-minimizing region R with $1 \geq |R| \geq \beta\delta > 0$ lying in the tangent cone at a vertex or the most complicated type of point that stays within a bounded distance of the origin. By Lemma 3.5, there are constants $C_1, \delta_1 > 0$ such that for α large, in any $\delta_1$-ball B, $S_\alpha = \partial R_\alpha$ minimizes area plus $C_1|V|$. By Corollary 3.7, the sum of the diameters of components of $S_\alpha$ and hence of the $R_\alpha$ is bounded. Since the polytope is being scaled up, by translation of components, we may assume that $R_\alpha$ is contained in the union of disjoint balls about the vertices. By Theorem 3.3.B, we may assume that $R_\alpha$ consists of balls about vertices v. Since for such a ball, $A = c_v V^t$ with $t = (n-2)/(n-1) < 1$, by the concavity of $V^t$, a single ball is best.

**3.9. Remark.** A more general but more technical version of the concentration argument at the beginning of the proof of Theorem 3.8 uses biLipschitz equivalence with $\mathbf{R}^n$ as in the proof of Lemma 3.6. The small polytopes $K_i$ may then be replaced by cubes in $\mathbf{R}^n$. This approach obviates the need for our cute Lemma 3.4.

**3.10. Nonconvex polytopes.** Theorem 3.8 fails for general "nonconvex polytopes." For example, consider two tall skinny pyramids in $\mathbf{R}^3$, attached at the apices. By Theorem 3.8, for small volume, in each pyramid separately a ball about the apex is perimeter minimizing. Hence in their union, those same sets are perimeter minimizing. But balls about the common apex are unions of such sets and have more perimeter.

For n = 2, Theorem 3.8 holds for nonconvex polyhedra locally homeomorphic to $\mathbf{R}^2$, for which the tangent cones are all isometric to circular cones. The only additional ingredient for the proof is that in a cone with cone angle greater than $2\pi$, Euclidean discs (away from the apex) are perimeter minimizing. This fact follows for example from Theorem 3.2.

For n ≥ 3, Theorem 3.8 fails even for nonconvex polytopes homeomorphic to the 3-sphere. For example, let $K^2$ be the surface of a unit cube in $\mathbf{R}^3 = \{x \in \mathbf{R}^4 : x_4 = 1\}$, surmounted by a tall skinny tetrahedron with apex p. Note that in the cone $0 \ast K^2$ over $K^2$, for small volume, a ball about a point q in $0 \ast \{p\}$ has less perimeter than a ball about 0, the apex of the cone, because of the small solid angle. To get the counterexample polytope $P^3$, truncate the cone at $\{x_4 = 1\}$ and surmount it on the surface of a big hypercube in $\mathbf{R}^4$. A ball about q still has less perimeter than a ball about any vertex of $P^3$, essentially because the exterior of K in $\mathbf{R}^3$ has no small solid angles.

**3.11. Solid polytopes.** The results of this paper apply to the nicer category of solid polytopes $P^n \subset \mathbf{R}^n$, which have no interior singularities. Theorem 3.3 with uniqueness and proof go over without smoothing to solid spherical polytopes $Q^{n-1} \subset \mathbf{S}^{n-1}$, a result apparently first proved (without uniqueness) by Lions and Pacella ([LP, Thm. 1.1], see also [M3, Rmk. after Thm. 10.6]). Likewise Theorem 3.8 with uniqueness and proof go over to solid polytopes $P^n \subset \mathbf{R}^n$.

Theorem 3.2 also has the following general consequence.

**3.12. Theorem.** *Let $M^n$ be a compact $C^1$ submanifold of $\mathbf{R}^{n+1}$ (with or without boundary), or polytope (convex or not, solid or boundary, but compact), or any compact Lipschitz neighborhood retract. Then in the cone over a small homothetic copy of M, balls about the apex are isoperimetric.*

*Proof.* The isoperimetric profile for M satisfies

$$I(M)(V) \geq cV^{(n-1)/n}$$

for $0 \leq V \leq |M|/2$ ([Fed, 4.4.2(2)]; see [M1, 12.3]). Hence the isoperimetric profile for a small homothetic copy M′ of M, renormalized with constant density to make $|M'| = |\mathbf{S}^n|$, satisfies

(1) $$I(M') \geq I(\mathbf{S}^n),$$

with equality only at the endpoints. Hence by Theorem 3.2, in the cone over M′, balls about the apex are uniquely perimeter minimizing.

*Remarks.* Similar results were already known for products M×**R** [DuS, 2.11] and apparently M×**R**$^n$ [Gon]. Such results follow from Theorem 3.2 for M×**R** and M×**S**$^1$ by comparison with results on **S**$^n$×**R** [P] and **S**$^n$×**S**$^1$ [PR].

V. Bayle pointed out to me that Berard, Besson, and Gallot [BBG, Thm. 2] give explicit sufficient conditions on the scaling down of M for (1) and hence Theorem 3.12 to hold, in terms of the Ricci curvature, the diameter, and the dimension of M.

**3.13. Conjecture** (M. Hutchings [H, Conj. 3]). *Let M, N be compact Riemannian manifolds. Then for fixed volume fraction, in the product of a small homothetic copy of M with N, an isoperimetric region is the product of M with an isoperimetric region in N.*